\documentclass[a4paper]{article}
\usepackage{amsmath,amssymb,mathrsfs,graphicx,psfrag}
\allowdisplaybreaks[1]

\DeclareMathOperator{\im}{Im}
\DeclareMathOperator{\fp}{f.p.}
\newtheorem{thm}{Theorem}
\begin{document}
\title{A numerical method for Hadamard finite-part integrals with an integral power singularity at an endpoint}
\author{Hidenori Ogata%
\footnote{
Department of Computer and Network Engineering, Graduate School of Informatics and Engineering, 
The University of Electro-Communications, 1-5-1 Chofugaoka, Chofu, Tokyo 182-8585, Japan, 
(e-mail) {\tt ogata@im.uec.ac.jp}
}}
\maketitle
\begin{abstract}
 In this paper, we propose a numerical method for computing Hadamard finite-part integrals with an integral-power 
 singularity at an endpoint, the part of the divergent integral which is finite as a limiting procedure. 
 In the proposed method, we express the desired finite-part integral using a complex loop integral, and obtain 
 the finite-part integral by evaluating the complex integral by the trapezoidal rule. 
 Theoretical error estimate and some numerical examples show the effectiveness of the proposed method. 
\end{abstract}
\section{Introduction}
\label{sec:introduction}
The integral 
\begin{equation*}
 \int_0^1 x^{-1}f(x)\mathrm{d}x, 
\end{equation*}
where $f(x)$ is an analytic function on the closed interval $[0,1]$ such that $f(0)\neq 0$, 
is divergent. 
However, if $f(x)$ is analytic on the closed interval $[0,1]$, 
for $\epsilon$ such that $0<\epsilon\ll 1$, we have the following using integral by part. 
\begin{align*}
 \int_{\epsilon}^{1}x^{-1}f(x)\mathrm{d}x = \: & 
 \int_{\epsilon}^{1}(\log x)^{\prime}f(x)\mathrm{d}x
 \\ 
 = \: & 
 \bigg[\log x f(x)\bigg]_{\epsilon}^{1} - \int_{\epsilon}^{1}\log x f^{\prime}(x)\mathrm{d}x
 \\ 
 = \: & 
 - f(\epsilon)\log\epsilon - \int_{\epsilon}^{1}\log x f^{\prime}(x)\mathrm{d}x
 \\ 
 = \: & - f(0)\log\epsilon + 
 \mbox{(terms finite as $\epsilon\downarrow 0$)}.
\end{align*}
Therefore, we can define the so-called Hadamard finite-part (f.p.) integral by 
\begin{equation*}
 \fp\int_0^1 x^{-1}f(x)\mathrm{d}x = 
  \lim_{\epsilon\downarrow 0}
  \left\{
   \int_{\epsilon}^1 x^{-1}f(x) + f(0)\log\epsilon
  \right\}
\end{equation*}
\cite{EstradaKanwal1989}. 
Similarly, we can define the f.p. integrals
\begin{equation}
 \label{eq:fp-integral0}
  \fp\int_0^1 x^{-n}f(x)\mathrm{d}x \quad ( \: n = 1, 2, \ldots \: ).
\end{equation}
We propose a numerical method for computing the f.p. integrals 
(\ref{eq:fp-integral0}). 
In the proposed method, we express the f.p. integral by a complex loop integral, and 
we obtain the f.p. integral by evaluating the complex integral by the trapezoidal formula with equal mesh.

Previous works related to this paper are as follows. Ogata and Hirayama proposed a numerical integration 
method based on hyperfunction theory, a theory of generalized functions based on complex function theory, 
where they obtain ordinary integrals by expressing them as complex integrals and evaluating them 
by the conventional numerical integral formulas \cite{OgataHirayama2018}. 
For Cauchy principal value integrals and Hadamard f.p. integrals with a singularity inside the integral interval 
\begin{equation}
 \label{eq:fp-integral01}
 \fp\int_0^1 \frac{f(x)}{(x-\lambda)^n}\mathrm{d}x \quad ( \: n = 1, 2, \ldots \: ),
\end{equation}
many approximation methods have been proposed. 
Elliot and Paget proposed a Gauss type numerical integration formula for Cauchy principal value integrals 
(\ref{eq:fp-integral01}) with $n=1$ \cite{ElliotPaget1979}, and 
Paget proposed a Gauss type formula for Hadamard finite-part integrals (\ref{eq:fp-integral01}) with $n=2$ 
\cite{Paget1981}. 
Bialecki proposed approximation formulas for (\ref{eq:fp-integral01}) based on the Sinc method 
\cite{Bialecki1990a,Bialecki1990b}, that is, 
methods using the trapezoidal formula together with variable transforms as in the DE formula \cite{TakahasiMori1974}. 
The author et al. improved them and 
proposed a DE-type numerical integration formula for Cauchy principal-value integrals and 
Hadamard finite-part integrals with an integral power singularity inside the integral interval 
\cite{OgataSugiharaMori2000}.

The remainder of this paper is structured as follows. 
In Section \ref{sec:fp-integral}, we define the f.p. integrals (\ref{eq:fp-integral0}) 
and show the expression of them by complex loop integral. 
Then, we give an approximation formula for the desired f.p. integral. 
In Section \ref{sec:example}, we show some numerical examples which show the effectiveness of the proposed method. 
In Section \ref{sec:summary}, we give a summary of this paper. 
\section{Hadamard finite-part integrals and its approximation}
\label{sec:fp-integral}
We define the Hadamard finite-part integrals (\ref{eq:fp-integral0}) by
\begin{multline}
 \label{eq:fp-integral}
 \fp\int_0^1 x^{-n}f(x)\mathrm{d}x 
 \\
 = 
 \lim_{\epsilon\downarrow 0}
 \left\{
 \int_{\epsilon}^1 x^{-n}f(x) 
 - \sum_{k=0}^{n-2}\frac{\epsilon^{k+1-n}}{k!(n-1-k)}f^{(k)}(0) 
 + \frac{\log\epsilon}{(n-1)!}f^{(n-1)}(0) 
 \right\}
 \\ 
 ( \: n = 1, 2, \ldots \: ),
\end{multline}
where the integrand $f(x)$ is analytic on the closed interval $[0,1]$, 
and the second term on the right-hand side is zero if $n=1$. 
We can show that it is well-defined using integral by part as follows. 
\begin{align*}
 & \int_{\epsilon}^{1}x^{-n}f(x)\mathrm{d}x 
 \\ 
 = \: & 
 - \frac{1}{n-1}\int_{\epsilon}^{1}(x^{-(n-1)})^{\prime}f(x)\mathrm{d}x
 \\ 
 = \: & 
 - \frac{1}{n-1}
 \left\{
 \bigg[x^{-(n-1)}f(x)\bigg]_{\epsilon}^{1} - 
 \int_{\epsilon}^{1}x^{-(n-1)}f^{\prime}(x)\mathrm{d}x
 \right\}
 \\ 
 = \: & 
 - \frac{f(1)}{n-1} + \frac{\epsilon^{1-n}}{n-1}f(\epsilon) + 
 \frac{1}{n-1}\int_{\epsilon}^{1}x^{-(n-1)}f^{\prime}(x)\mathrm{d}x
 \\ 
 = \: & 
 \frac{\epsilon^{1-n}}{n-1}\sum_{k=0}^{n-2}
 \frac{\epsilon^k}{k!}f^{(k)}(0) 
 - \frac{1}{(n-1)(n-2)}\int_{\epsilon}^{1}
 (x^{-(n-2)})^{\prime}f^{\prime}(x)\mathrm{d}x
 \\ 
 & + \mbox{(terms finite as $\epsilon\downarrow 0$, which is denoted by ^^ ^^ $\cdots$'' below)}
 \\ 
 = \: & 
 \frac{\epsilon^{1-n}}{n-1}\sum_{k=0}^{n-2}
 \frac{\epsilon^k}{k!}f^{(k)}(0)
 +
 \frac{\epsilon^{2-n}}{(n-1)(n-2)}\sum_{k=0}^{n-3}
 \frac{\epsilon^k}{k!}f^{(k+1)}(0)
 \\  
 & 
 - 
 \frac{1}{(n-1)(n-2)(n-3)}\int_{\epsilon}^{1}
 (x^{-(n-3)})^{\prime}f^{\prime\prime}(x)\mathrm{d}x
 + \cdots 
 \\ 
 = \: & 
 \cdots
 \\ 
 = \: & 
 \frac{\epsilon^{1-n}}{n-1}\sum_{k=0}^{n-2}
 \frac{\epsilon^k}{k!}f^{(k)}(0)
 +
 \frac{\epsilon^{2-n}}{(n-1)(n-2)}\sum_{k=0}^{n-3}
 \frac{\epsilon^k}{k!}f^{(k+1)}(0)
 \\ 
 & + 
 \frac{\epsilon^{3-n}}{(n-1)(n-2)(n-3)}\sum_{k=0}^{n-4}
 \frac{\epsilon^k}{k!}f^{(k+2)}(0) + \cdots + \frac{\epsilon^{-1}}{(n-1)!}f^{(n-2)}(0)
 - \frac{\log\epsilon}{(n-1)!}f^{(n-1)}(0)  
 \\ 
 & 
 - \frac{1}{(n-1)!}\int_{\epsilon}^{1}\log x f^{(n)}(x)\mathrm{d}x + \cdots 
 \\ 
 = \: & 
 \frac{\epsilon^{1-n}}{n-1}\sum_{k=0}^{n-2}
 \frac{\epsilon^k}{k!}f^{(k)}(0)
 +
 \frac{\epsilon^{2-n}}{(n-1)(n-2)}\sum_{k=1}^{n-2}
 \frac{\epsilon^{k-1}}{(k-1)!}f^{(k)}(0)
 \\ 
 & + 
 \frac{\epsilon^{3-n}}{(n-1)(n-2)(n-3)}\sum_{k=2}^{n-2}
 \frac{\epsilon^{k-2}}{(k-2)!}f^{(k)}(0) + \cdots + \frac{\epsilon^{-1}}{(n-1)!}f^{(n-2)}(0)
 - \frac{\log\epsilon}{(n-1)!}f^{(n-1)}(0)  
 \\ 
 & 
 - \frac{1}{(n-1)!}\int_{\epsilon}^{1}\log x f^{(n)}(x)\mathrm{d}x + \cdots 
 \\ 
 = \: & 
 \frac{\epsilon^{1-n}}{n-1}f(0) 
 + 
 \epsilon^{2-n}f^{\prime}(0)
 \left\{\frac{1}{n-1} + \frac{1}{(n-1)(n-2)}\right\}
 \\ 
 & 
 + 
 \frac{\epsilon^{1-n}}{n-1}\sum_{k=2}^{n-2}\frac{\epsilon^k}{k!}f^{(k)}(0)
 +
 \frac{\epsilon^{2-n}}{(n-1)(n-2)}\sum_{k=2}^{n-2}
 \frac{\epsilon^{k-1}}{(k-1)!}f^{(k)}(0)
 \\ 
 & + 
 \frac{\epsilon^{3-n}}{(n-1)(n-2)(n-3)}\sum_{k=2}^{n-2}
 \frac{\epsilon^{k-2}}{(k-2)!}f^{(k)}(0) + \cdots + \frac{\epsilon^{-1}}{(n-1)!}f^{(n-2)}(0)
 - \frac{\log\epsilon}{(n-1)!}f^{(n-1)}(0)  
 \\ 
 & 
 - \frac{1}{(n-1)!}\int_{\epsilon}^{1}\log x f^{(n)}(x)\mathrm{d}x + \cdots 
 \\ 
 = \: & 
 \frac{\epsilon^{1-n}}{n-1}f(0) + \frac{\epsilon^{2-n}}{n-2}f^{\prime}(0)
 + 
 \frac{\epsilon^{3-n}}{n-1}f^{\prime\prime}(0)
 \left\{
 \frac{1}{2!} + \frac{1}{n-2} + \frac{1}{(n-2)(n-3)}
 \right\}
 \\ 
 & 
 + 
 \frac{\epsilon^{1-n}}{n-1}\sum_{k=3}^{n-2}\frac{\epsilon^k}{k!}f^{(k)}(0)
 +
 \frac{\epsilon^{2-n}}{(n-1)(n-2)}\sum_{k=3}^{n-2}
 \frac{\epsilon^{k-1}}{(k-1)!}f^{(k)}(0)
 \\ 
 & + 
 \frac{\epsilon^{3-n}}{(n-1)(n-2)(n-3)}\sum_{k=3}^{n-2}
 \frac{\epsilon^{k-2}}{(k-2)!}f^{(k)}(0) + \cdots + \frac{\epsilon^{-1}}{(n-1)!}f^{(n-2)}(0)
 - \frac{\log\epsilon}{(n-1)!}f^{(n-1)}(0)  
 \\ 
 & + \cdots 
 \\ 
 = \: & 
 \frac{\epsilon^{1-n}}{n-1}f(0) + \frac{\epsilon^{2-n}}{n-2}f^{\prime}(0)
 + 
 \frac{\epsilon^{3-n}}{2!(n-3)}f^{\prime\prime}(0)
 \\ 
 & 
 + 
 \frac{\epsilon^{4-n}}{n-1}f^{\prime\prime\prime}(0)
 \bigg\{
 \frac{1}{3!} + 
 \underbrace{\frac{1}{2!(n-2)} + 
 \underbrace{\frac{1}{(n-2)(n-3)} + \frac{1}{(n-2)(n-3)(n-4)}}_{\frac{1}{(n-2)(n-4)}}
 }_{\frac{1}{2!(n-4)}}
 \bigg\}
 \\ 
 & 
 + 
 \frac{\epsilon^{1-n}}{n-1}\sum_{k=4}^{n-2}\frac{\epsilon^k}{k!}f^{(k)}(0)
 +
 \frac{\epsilon^{2-n}}{(n-1)(n-2)}\sum_{k=4}^{n-2}
 \frac{\epsilon^{k-1}}{(k-1)!}f^{(k)}(0)
 \\ 
 & + 
 \frac{\epsilon^{3-n}}{(n-1)(n-2)(n-3)}\sum_{k=4}^{n-2}
 \frac{\epsilon^{k-2}}{(k-2)!}f^{(k)}(0) 
 \\ 
 & 
 +
 \frac{\epsilon^{4-n}}{(n-1)(n-2)(n-3)(n-4)}\sum_{k=4}^{n-2}
 \frac{\epsilon^{k-3}}{(k-3)!}f^{(k)}(0)
 +
 \cdots + \frac{\epsilon^{-1}}{(n-1)!}f^{(n-2)}(0)
 \\ 
 & 
 - \frac{\log\epsilon}{(n-1)!}f^{(n-1)}(0)   
 + \cdots 
 \\ 
 = \: & 
 \frac{\epsilon^{1-n}}{n-1}f(0) + \frac{\epsilon^{2-n}}{n-2}f^{\prime}(0)
 + 
 \frac{\epsilon^{3-n}}{2!(n-3)}f^{\prime\prime}(0)
 + 
 \frac{\epsilon^{4-n}}{3!(n-4)}f^{\prime\prime\prime}(0)
 \\ 
 & 
 + 
 \frac{\epsilon^{1-n}}{n-1}\sum_{k=4}^{n-2}\frac{\epsilon^k}{k!}f^{(k)}(0)
 +
 \frac{\epsilon^{2-n}}{(n-1)(n-2)}\sum_{k=4}^{n-2}
 \frac{\epsilon^{k-1}}{(k-1)!}f^{(k)}(0)
 \\ 
 & + 
 \frac{\epsilon^{3-n}}{(n-1)(n-2)(n-3)}\sum_{k=4}^{n-2}
 \frac{\epsilon^{k-2}}{(k-2)!}f^{(k)}(0) 
 \\ 
 & 
 +
 \frac{\epsilon^{4-n}}{(n-1)(n-2)(n-3)(n-4)}\sum_{k=4}^{n-2}
 \frac{\epsilon^{k-3}}{(k-3)!}f^{(k)}(0)
 +
 \cdots + \frac{\epsilon^{-1}}{(n-1)!}f^{(n-2)}(0)
 \\ 
 & 
 - \frac{\log\epsilon}{(n-1)!}f^{(n-1)}(0)   
 + \cdots 
 \\ 
 = \: & 
 \frac{\epsilon^{1-n}}{n-1}f(0) + \frac{\epsilon^{2-n}}{n-2}f^{\prime}(0)
 + 
 \frac{\epsilon^{3-n}}{2!(n-3)}f^{\prime\prime}(0)
 + 
 \frac{\epsilon^{4-n}}{3!(n-4)}f^{\prime\prime\prime}(0)
 + \cdots + 
 \frac{\epsilon^{-1}}{(n-2)!}f^{(n-2)}(0) 
 \\ 
 & - \frac{\log\epsilon}{(n-1)!}f^{(n-1)}(0) + \cdots.
\end{align*}

\begin{figure}[htbp]
 \begin{center}
  \psfrag{0}{$\mathrm{O}$}
  \psfrag{1}{$1$}
  \psfrag{p}{$\Gamma_{\epsilon}^{(+)}$}
  \psfrag{m}{$\Gamma_{\epsilon}^{(-)}$}
  \psfrag{a}{$C_{\epsilon}^{(0)}$}
  \psfrag{b}{$C_{\epsilon}^{(1)}$}
  \psfrag{e}{$\epsilon$}
  \psfrag{f}{$1-\epsilon$}
  \includegraphics[width=0.7\textwidth]{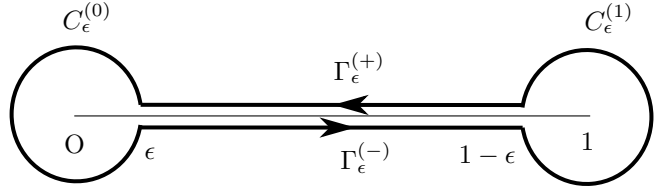}
 \end{center}
 \caption{The integral paths.}
 \label{fig:integral-path}
\end{figure}

The f.p. integral is expressed using a complex loop integral as in the following theorem. 
\begin{thm}
 \label{thm:complex-integral}
 We suppose that $f(z)$ is analytic in a complex domain $D$ containing the closed interval $[0,1]$ 
 in its interior. Then, the f.p. integral (\ref{eq:fp-integral}) is expressed as 
 \begin{multline}
  \label{eq:complex-integral}
   \fp\int_0^1 x^{-n}f(x)\mathrm{d}x
   = 
   \frac{1}{2\pi\mathrm{i}}\oint_C 
   z^{-n}f(z)\log\left(\frac{z}{z-1}\right)\mathrm{d}z 
   - 
   \sum_{k=0}^{n-2}\frac{f^{(k)}(0)}{k!(n-1-k)} \\ ( \: n = 1, 2, \ldots \: ),
 \end{multline}
 where $C$ is a closed complex integral path in $D$ encircling the interval $[0,1]$ in the positive sense, 
 and the second term on the right-hand side of (\ref{eq:complex-integral}) is zero if $n=1$. 
\end{thm}
\paragraph{Proof of Theorem \ref{thm:complex-integral}}
Using Cauchy's integral theorem, we have
\begin{multline}
 \label{eq:proof-complex-integral}
 \frac{1}{2\pi\mathrm{i}}\int_C
 z^{-n}f(z)\log\left(\frac{z}{z-1}\right)\mathrm{d}z
 \\
 = 
 \frac{1}{2\pi\mathrm{i}}
 \left(
 \int_{C_{\epsilon}^{(0)}} + \int_{C_{\epsilon}^{(1)}} + 
 \int_{\Gamma_{\epsilon}^{(+)}} + \int_{\Gamma_{\epsilon}^{(-)}}
 \right)
 z^{-n}f(z)\log\left(\frac{z}{z-1}\right)\mathrm{d}z, 
\end{multline}
where the integral paths $C_{\epsilon}^{(0)}$, $C_{\epsilon}^{(1)}$, 
$\Gamma_{\epsilon}^{(+)}$ and $\Gamma_{\epsilon}^{(-)}$ are respectively 
\begin{align*}
 C_{\epsilon}^{(0)} = \: & 
 \{ \: \epsilon\mathrm{e}^{\mathrm{i}\theta} \: | \: 0 \leqq \theta \leqq \pi \: \}, 
 \\ 
 C_{\epsilon}^{(1)} = \: & 
 \{ \: 1 + \epsilon\mathrm{e}^{\mathrm{i}\theta} \: | \: 0 \leqq \theta \leqq \pi \: \}, 
 \\ 
 \Gamma_{\epsilon}^{(+)} = \: & 
 \{ \: x \in\mathbb{R} \: | \: 1-\epsilon \geqq x \geqq \epsilon \: \}, 
 \\ 
 \Gamma_{\epsilon}^{(-)} = \: & 
 \{ \: x \in\mathbb{R} \: | \: \epsilon \leqq x \leqq 1-\epsilon \: \}
\end{align*}
with small $\epsilon>0$ (see Figure \ref{fig:integral-path}).
\begin{table}[htbp]
 \caption{The arguments of the functions appearing in the complex integral (\ref{eq:proof-complex-integral}).}
 \begin{center}
  \begin{tabular}{c|ccccc}
   \hline 
   $z$ & $\Gamma_{\epsilon}^{(+)}$ & 
	   $C_{\epsilon}^{(0)}$ & 
	       $\Gamma_{\epsilon}^{(-)}$ & 
		   $C_{\epsilon}^{(1)}$ & 
		       $\Gamma_{\epsilon}^{(+)}$ \\
   \hline
   $\arg z$ & $0$ & & $2\pi$ & & $2\pi$ \\ 
   $\arg(z-1)$ & $\pi$ & & $\pi$ & & $3\pi$ \\ 
   $\arg(z/(z-1))$ & $-\pi$ & & $\pi$ & & $-\pi$ \\ 
   \hline
  \end{tabular}
 \end{center}
 \label{tab:argument}
\end{table}
From Table \ref{tab:argument}, we have
\begin{align*}
 & 
 \frac{1}{2\pi\mathrm{i}}
 \left( \int_{\Gamma_{\epsilon}^{(+)}} + \int_{\Gamma_{\epsilon}^{(-)}}\right)
 z^{-n}f(z)\log\left(\frac{z}{z-1}\right)\mathrm{d}z
 \\ 
 = \: & 
 - \frac{1}{2\pi\mathrm{i}}\int_{\epsilon}^{1-\epsilon}
 x^{-n}f(x)\left\{\log\left(\frac{x}{1-x} - \mathrm{i}\pi\right)\right\}\mathrm{d}x
 + 
 \frac{1}{2\pi\mathrm{i}}\int_{\epsilon}^{1-\epsilon}
 x^{-n}f(x)\left\{\log\left(\frac{x}{1-x} + \mathrm{i}\pi\right)\right\}\mathrm{d}x
 \\ 
 = \: & 
 \int_{\epsilon}^{1-\epsilon}x^{-n}f(x)\mathrm{d}x 
 = 
 \int_{\epsilon}^{1}x^{-n}f(x)\mathrm{d}x + \mathrm{O}(\epsilon).
\end{align*}
The integral on $C_{\epsilon}^{(0)}$ is written as 
\begin{align*}
 & 
 \frac{1}{2\pi\mathrm{i}}\int_{C_{\epsilon}^{(0)}}
 z^{-n}f(z)\log\left(\frac{z}{z-1}\right)\mathrm{d}z
 \\ 
 = \: & 
 \frac{1}{2\pi\mathrm{i}}\int_0^{2\pi}
 \epsilon^{-n}\mathrm{e}^{-\mathrm{i}n\theta}f(\epsilon\mathrm{e}^{\mathrm{i}\theta})
 \log\left(\frac{-\epsilon\mathrm{e}^{\mathrm{i}\theta}}{1-\epsilon\mathrm{e}^{\mathrm{i}\theta}}\right)
 \mathrm{i}\epsilon\mathrm{e}^{\mathrm{i}\theta}\mathrm{d}\theta
 \\ 
 = \: & 
 \frac{\epsilon^{1-n}}{2\pi}\int_0^{2\pi}f(\epsilon\mathrm{e}^{\mathrm{i}\theta})
 \log(\epsilon\mathrm{e}^{\mathrm{i}(\theta-\pi)})\mathrm{e}^{-\mathrm{i}(n-1)\theta}
 \mathrm{d}\theta
 - 
 \frac{\epsilon^{1-n}}{2\pi}\int_0^{2\pi}f(\epsilon\mathrm{e}^{\mathrm{i}\theta})
 \log(1-\epsilon\mathrm{e}^{\mathrm{i}\theta})\mathrm{e}^{-\mathrm{i}(n-1)\theta}
 \mathrm{d}\theta.
\end{align*}
The first integral on the right-hand side is written as
\begin{align*}
 & 
 \frac{\epsilon^{1-n}}{2\pi}\int_0^{2\pi}
 f(\epsilon\mathrm{e}^{\mathrm{i}\theta})\log(\epsilon\mathrm{e}^{\mathrm{i}(\theta-\pi)})
 \mathrm{e}^{-\mathrm{i}(n-1)\theta}
 \mathrm{d}\theta
 \\ 
 = \: & 
 \frac{\epsilon^{1-n}}{2\pi}\int_0^{2\pi}
 \left\{\sum_{k=0}^{\infty}\frac{\epsilon^k}{k!}f^{(k)}(0)\mathrm{e}^{\mathrm{i}k\theta}\right\}
 \left\{\log\epsilon + \mathrm{i}(\theta-\pi)\right\}
 \mathrm{e}^{-\mathrm{i}(n-1)\theta}
 \mathrm{d}\theta
 \\ 
 = \: & 
 \frac{\epsilon^{1-n}}{2\pi}\log\epsilon\sum_{k=0}^{\infty}
 \int_0^{2\pi}\mathrm{e}^{\mathrm{i}(k-n+1)\theta}
 \mathrm{d}\theta
 + 
 \frac{\mathrm{i}(n-1)\epsilon^{1-n}}{2\pi}\sum_{k=0}^{\infty}
 \frac{\epsilon^k}{k!}f^{(k)}(0)
 \int_0^{2\pi}(\theta-\pi)\mathrm{e}^{\mathrm{i}(k-n+1)\theta}
 \mathrm{d}\theta
 \\ 
 = \: & 
 \frac{\log\epsilon}{(n-1)!}f^{(n-1)}(0)
 - 
 \sum_{k=0}^{n-2}\frac{\epsilon^{k-n+1}}{k!(n-1-k)}f^{(k)}(0) + \mathrm{O}(\epsilon), 
\end{align*}
where we exchanged the order of the integral and the infinite summation on the second equality 
since the infinite sum is uniformly convergent on $0\leqq\theta\leqq 2\pi$. 
Similarly, the second integral is written as
\begin{equation*}
 \frac{\epsilon^{1-n}}{2\pi}\int_0^{2\pi}
 f(\epsilon\mathrm{e}^{\mathrm{i}\theta})\log(1-\epsilon\mathrm{e}^{\mathrm{i}\theta})
 \mathrm{e}^{-\mathrm{i}(n-1)\theta}
 \mathrm{d}\theta
 = 
 - \sum_{k=0}^{n-2}\frac{f^{(k)}(0)}{k!(n-1-k)}.
\end{equation*}
Then, we have
\begin{multline*}
 \frac{1}{2\pi\mathrm{i}}\int_{C_{\epsilon}^{(0)}}
 z^{-n}f(z)\log\left(\frac{z}{z-1}\right)
 \mathrm{d}z
 \\ 
 = 
 - \sum_{k=0}^{n-2}\frac{\epsilon^{k-n+1}}{k!(n-1-k)}f^{(k)}(0) 
 + \frac{\log\epsilon}{(n-1)!}f^{(n-1)}(0)
 + \sum_{k=0}^{n-2}\frac{f^{(k)}(0)}{k!(n-1-k)}
 + \mathrm{O}(\epsilon).
\end{multline*}
As to the integral on $C_{\epsilon}^{(1)}$, we have
\begin{equation*}
 \frac{1}{2\pi\mathrm{i}}\int_{C_{\epsilon}^{(1)}}
  z^{-n}f(z)\log\left(\frac{z}{z-1}\right)\mathrm{d}z
  = 
  \mathrm{O}(\epsilon\log\epsilon)
\end{equation*}
since $z^{-n}f(z)\log z$ is analytic near $z=1$. 
Summarizing the above calculations, we have
\begin{multline*}
 \frac{1}{2\pi\mathrm{i}}\oint_C
 z^{-n}f(z)\log\left(\frac{z}{z-1}\right)
 \mathrm{d}z
 \\ 
 = 
 \int_{\epsilon}^{1}x^{-n}f(x)\mathrm{d}x
 - 
 \sum_{k=0}^{n-2}\frac{\epsilon^{k-n+1}}{k!(n-1+k)}f^{(k)}(0) 
 \\ 
 + 
 \frac{\log\epsilon}{(n-1)!}f^{(n-1)}(0)
 + 
 \sum_{k=0}^{n-2}\frac{f^{(k)}(0)}{k!(n-1+k)}
 + 
 \mathrm{O}(\epsilon\log\epsilon).
\end{multline*}
Taking the limit $k\downarrow 0$, we have (\ref{eq:complex-integral}). 

\hfill\rule{1.5ex}{1.5ex}

The complex integral in (\ref{eq:complex-integral}) is the integral of an analytic function 
over an interval of the length of one period, and it is accurately approximated by the trapezoidal 
formula with equal mesh. 
Using a parameterization of the closed integral path 
\begin{equation*}
 C \: : \: z = \varphi(u), \quad 0 \leqq u \leqq u_{\rm p}, 
\end{equation*}
where $\varphi(u)$ is a periodic function of period $u_{\rm p}$, we obtain the following approximation 
formula for the f.p. integral.
\begin{multline}
 \label{eq:approx-fp-integral}
 \fp\int_0^1 x^{-n}f(x)\mathrm{d}x 
 \simeq 
 I_N^{(n)}[f]
 \\ 
 \equiv 
 \frac{h}{2\pi\mathrm{i}}\sum_{k=0}^{N-1}\varphi(kh)^{-n}f(\varphi(kh))
 \log\left(\frac{\varphi(kh)}{\varphi(kh)-1}\right)\varphi^{\prime}(kh)
 - 
 \sum_{k=0}^{n-2}\frac{f^{(k)}(0)}{k!(n-1-k)}
 \\ 
 \left( \: h = \frac{u_{\rm p}}{N} \: \right), 
\end{multline}
where the second term on the right-hand side is zero if $n=1$. 

We remark here that, if the integrand $f(x)$ is real valued on $[0,1]$ and the integral path $C$ is symmetric 
with respect to the real axis, we can reduce the number of sampling points $N$ by half. 
In fact, in this case, we have 
$f(\overline{z}) = \overline{f(z)}$ due to the reflection principle, and 
$\varphi(-u) = \overline{\varphi(u)}$, $\varphi^{\prime}(-u) = - \overline{\varphi^{\prime}(u)}$. 
Then, we have
\begin{align}
 \nonumber
 & 
 \fp\int_0^1 x^{-n}f(x)\mathrm{d}x 
 \simeq 
 {I^{\prime}}_{N}^{(n)}[f]
 \\ 
 \nonumber
 \equiv \: & 
 \frac{h}{2\pi}\im
 \left\{
 \varphi(0)^{-n}f(\varphi(0))
 \log\left(\frac{\varphi(0)}{1-\varphi(0)}\right)\varphi^{\prime}(0)
 \right.
 \\
 \nonumber
 & 
 \hspace{12mm}
 \left.
 +
 \varphi\left(\frac{u_{\rm p}}{2}\right)^{-n}
 f\left(\varphi\left(\frac{u_{\rm p}}{2}\right)\right)
 \log\left(\frac{\varphi(u_{\rm p}/2)}{1-\varphi(u_{\rm p}/2)}\right)
 \varphi^{\prime}\left(\frac{u_{\rm p}}{2}\right)
 \right\}
 \\
 \nonumber
 & 
 +
 \frac{h}{\pi}\im\left\{\sum_{k=1}^{N-1}\varphi(kh)^{-n}f(\varphi(kh))
 \log\left(\frac{\varphi(kh)}{\varphi(kh)-1}\right)\varphi^{\prime}(kh)
 \right\}
 \\
 \label{eq:approx-fp-integral2}
 & 
 - 
 \sum_{k=0}^{n-2}\frac{f^{(k)}(0)}{k!(n-1-k)}
 \quad 
 \left( \: h = \frac{u_{\rm p}}{2N} \: \right), 
\end{align}

Applying the theorem in \S 4.6.5 in \cite{DavisRabinowitz1984} to the approximation of 
the complex integral by the trapezoidal formula in (\ref{eq:approx-fp-integral}), 
we have the following theorem 
on the error estimate of the approximation formula (\ref{eq:approx-fp-integral}).
\begin{thm}
 We suppose that 
 \begin{itemize}
  \item the strip domain 
	\begin{equation*}
	 D_d = 
	  \{ \: w \in\mathbb{C} \: | \: |\im w| < d \: \} 
	  \quad ( \: d > 0 \: )
	\end{equation*}
	is contained in $\mathbb{C}\setminus [0,1]$, 
  \item the parameterization function $\varphi(w)$ of $C$ is analytic in $D_d$, and 
  \item the integrand $f(z)$ is analytic in 
	\begin{equation*}
	 \varphi(D_d) = 
	  \{ \: \varphi(w) \: | \: w \in D_d \: \}.
	\end{equation*}
 \end{itemize}
 Then, we have the following inequality for arbitrary $0 < d^{\prime} < d$. 
 \begin{gather}
  \left|\fp\int_0^1 x^{-n}f(x)\mathrm{d}x - I_N^{(n)}[f] \right| 
  \leqq 
  \frac{d}{\pi}\mathscr{N}(f,n,d^{\prime})
  \frac{\exp(-2\pi d^{\prime}N/u_{\rm p})}{1-\exp(-2\pi d^{\prime}N/u_{\rm p})}, 
  \intertext{where}
  \mathscr{N}(f,n,d^{\prime}) = 
  \max_{|\im w|=d^{\prime}}
  \left|
  \varphi(w)^{-n}f(\varphi(w))\log\left(\frac{\varphi(w)}{1-\varphi(w)}\right).
  \right|
 \end{gather}
\end{thm}
This theorem says that the approximation (\ref{eq:approx-fp-integral}) converges exponentially 
as $N$ increases if the integrand function $f(x)$ is analytic on $[0,1]$ and the integral path $C$ is 
an analytic curve. 
%%%%%%%%%%%%%%%%%%%%%%%%%%%%%%%%%%%%%%%%%%%%%%%%%%%%%%%%%%%%%%%%%%%%%%%%%%%%%%%%%%%%%%%%%%%%%%%%%%%%%
\section{Numerical examples}
\label{sec:example}
\begin{figure}[htbp]
 \begin{center}
  \begin{tabular}{cc}
   \includegraphics[width=0.45\textwidth]{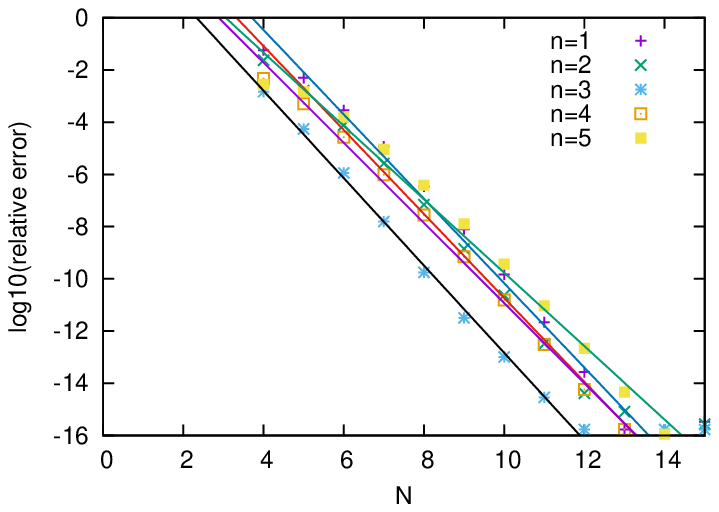} & 
       \includegraphics[width=0.45\textwidth]{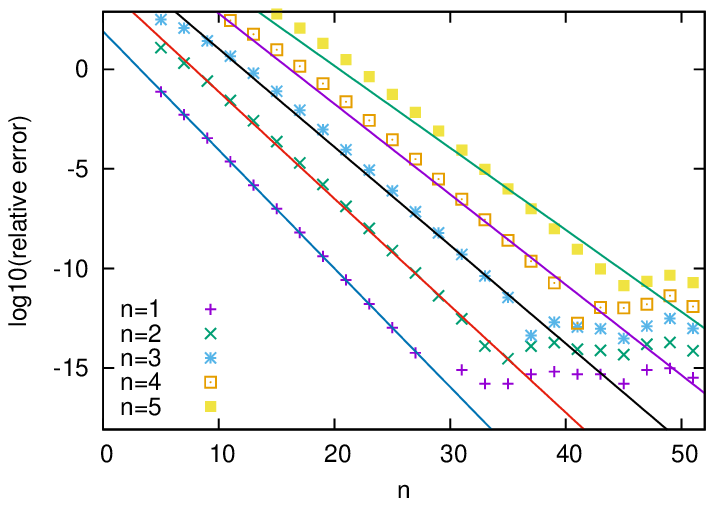}
       \\
   (1) & (2)
  \end{tabular}
 \end{center}
 \caption{The relative errors of the proposed method for f.p. integrals applied to the integrals (1) and (2).}
 \label{fig:example}
\end{figure}

We computed the integrals 
\begin{align*}
 \mathrm{(1)} \quad & 
 \fp\int_0^1 x^{-n}\mathrm{e}^x\mathrm{d}x = 
 \sum_{k=0 (k\neq n-1)}^{\infty}\frac{1}{k!(k-n+1)}, 
 \\ 
 \mathrm{(2)} \quad & 
 \fp\int_0^1 \frac{x^{-n}}{1+x}\mathrm{d}x = 
 (-1)^n\left\{\log 2 + \sum_{l=1}^{n-1}\frac{(-1)^l}{l}\right\},
\end{align*}
where the second term on the right-hand side of the integral (2) is zero if $n=1$, 
by the proposed method. 
All the computations were performed using programs coded in C++ with double precision working. 
The complex integral path $C$ was taken as the ellipse 
\begin{equation*}
 C \: : \: z = 
  \frac{1}{2} + \frac{1}{4}\left(\rho+\frac{1}{\rho}\right)\cos u + 
  \frac{1}{4}\left(\rho-\frac{1}{\rho}\right)\sin u, \quad 
  0 \leqq u \leqq 2\pi \quad ( \: \rho > 1 \: ), 
\end{equation*}
where the parameter $\rho$ is taken as $\rho=10$ for the integral (1) and 
$\rho=2$ for the integral (2). 
Figure \ref{fig:example} shows the relative errors of the approximation formula (\ref{eq:approx-fp-integral2}) 
applied to the integrals (1) and (2) as functions of the number of sampling points $N$. 
From these figures, the errors decay exponentially as $N$ increases. 
Table \ref{tab:example} shows the decay rates of the errors of the proposed method. 
\begin{table}[htbp]
 \caption{The decay rates of the errors of the proposed method for f.p. integrals applied to the integrals (1) and (2).}
 \begin{center}
  \begin{tabular}{ccccc}
   \hline
   $n$ &  & 1 & 2 & 3 \\
   \hline
   relative error & integral (1) & $\mathrm{O}(0.024^N)$ & 
	       $\mathrm{O}(0.025^N)$ & $\mathrm{O}(0.021^N)$ 
		   \\
   \cline{2-5}
   & integral (2) & $\mathrm{O}(0.25^N)$ & $\mathrm{O}(0.29^N)$ & $\mathrm{O}(0.32^N)$
		   \\
   \hline 
   &  & 4 & 5 \\
   \cline{3-5}
   &  & $\mathrm{O}(0.029^N)$ & $\mathrm{O}(0.039^N)$ 
	       \\
   \cline{3-5}
   &  & $\mathrm{O}(0.35^N)$ & $\mathrm{O}(0.38^N)$ 
	       \\
   \cline{3-5} 
  \end{tabular}
 \end{center}
 \label{tab:example}
\end{table}
%%%%%%%%%%%%%%%%%%%%%%%%%%%%%%%%%%%%%%%%%%%%%%%%%%%%%%%%%%%%%%%%%%%%%%%%%%%%%%%%%%%%%%%%%%%%%%%%%%%%%
\section{Summary}
\label{sec:summary}
We proposed a numerical method for Hadamard finite-part integrals with an integral order power singularity 
at an endpoint over a finite interval.  
In the proposed method, we express the desired f.p. integral using a complex loop integral, 
and we obtain the f.p. integral by evaluating the complex integral by the trapezoidal formula. 
Theoretical error estimate and numerical examples show the exponential convergence of the proposed method. 

We can also give approximation methods for f.p. integrals with a non-integral power singularity 
and f.p. integrals over a half-infinite interval in a way similar to this paper. 
They will be reported in other papers. 
%%%%%%%%%%%%%%%%%%%%%%%%%%%%%%%%%%%%%%%%%%%%%%%%%%%%%%%%%%%%%%%%%%%%%%%%%%%%%%%%%%%%%%%%%%%%%%%%%%%%%
\bibliographystyle{plain}
\bibliography{arxiv2019_2}
\end{document}